\documentclass[10pt]{article}
\usepackage{graphicx}
\usepackage{amsmath,amssymb,amsthm,amsfonts}
\usepackage{amssymb}
\newtheorem{thm}{Theorem}[section]

\newtheorem{lem}[thm]{Lemma}

\theoremstyle{definition}

\theoremstyle{remark}
\newtheorem{rem}{Remark}[section]

\numberwithin{equation}{section}

\DeclareMathSymbol{\C}{\mathalpha}{AMSb}{"43}

\newcommand{\eps}{\varepsilon}

\newcommand{\R}{{\mathbb{R}}}

\newcommand{\inte}{\int_{\mathbb{R}^2}}

\newcommand{\bsub}{\begin{subequations}}
\newcommand{\esub}{\end{subequations}$\!$}

\begin{document}

\title{Concentration behavior of standing waves for almost mass critical nonlinear Schr\"{o}dinger equations {\thanks{Email: YJG:yjguo@wipm.ac.cn; XYZ:zeng6656078@126.com; HSZ:hszhou@wipm.ac.cn; Tel:+86-2787199196}}}

\author{\small Yujin Guo, Xiaoyu Zeng  and Huan-Song Zhou \\
 \ \\
\small \it Wuhan Institute of Physics and Mathematics,
    Chinese Academy of Sciences,\\
    \small \it P.O. Box 71010, Wuhan 430071, P. R. China}

\date{}

\smallbreak \maketitle

\begin{abstract} We study  the following nonlinear Schr\"{o}dinger equation
\[
iu_t=-\Delta u+V(x)u-a|u|^qu \quad (t,x)\in \R^1\times \R^2,
\]
where $a>0, \ q\in(0,2)$ and $V(x)$ is some type of trapping potentials.  For any fixed $a>a^*:= \|Q\|_2^2$, where $Q$ is the unique
(up to translations) positive radial solution of $\Delta u-u+u^3=0$
in $\R^2$, by directly using constrained variational method and energy estimates we present a detailed analysis of the concentration and symmetry breaking  of the standing waves for the above equation as  $q\nearrow 2$.
\end{abstract}

\vskip 0.2truein



\vskip 0.2truein

\section{Introduction}

In this paper, we study the concentration and symmetry breaking of
standing waves for the following nonlinear Schr\"{o}dinger equation (NLS)
with a trapping potential and an attractive nonlinearity
\begin{equation}\label{new1}
iu_t=-\Delta u+V(x)u-a|u|^qu \quad\ (t,x)\in \R^1\times \R^2,
\end{equation}
where $a>0$, $0<q<2$, and $V(x)$ is a trapping potential.
The equation (\ref{new1}) with  $q=2$ arises in
Bose-Einstein condensates (BEC) as well as nonlinear optics, which has
been studied widely in recent years, see for examples, \cite{ca,GS,LSY,M,S}.
 In fact, when $q=2$ the above equation (\ref{new1}) is  the so-called mass critical NLS in $\R^2$, so $q=2$ is usually called a mass critical exponent for (\ref{new1}). Our this paper is focussed on the case where $q$ approaches $2$ from the left ($q\nearrow 2$, in short), which is what we mean by the almost mass critical NLS.

For (\ref{new1}), the standing
waves are the solutions of (\ref{new1}) with the form:
$u(t,x)=e^{i\omega t}\varphi_\omega (x)$, which implies that
$\varphi_\omega (x)$ satisfies the following elliptic partial
differential equation
\begin{equation}\label{1:pde2}
-\Delta u +(V+\omega )u -a|u
|^qu =0 \quad \text{in}\quad  \R^2.
\end{equation}
When $q=2$, (\ref{1:pde2}) is also called the time-independent Gross-Pitaevskii (GP) equation of Bose-Einstein condensates, where $\omega$ represents the chemical potential,  $V$ is an external potential, and $a$ is a coupling constant related to the number of bosons in a quantum system. Here $a>0 (resp. <0)$ means that the BEC is attractive ($resp.$ repulsive). In this paper, we consider only the attractive case, $i.e.,$ $a>0$.
It is well known that a minimizer of the following Gross-Pitaevskii (GP) energy functional
\begin{equation}\label{E}
E_q (u ):=\int_{\mathbb{R}^2}(|\nabla
u(x )|^2+V(x)|u(x )|^2)dx-\frac{2a}{q+2}\int_{\mathbb{R}^2}|u(x )|^{q+2}dx
\end{equation}
under the following
constraint
\begin{equation}\label{1:unit}
\inte u^2dx=1
\end{equation}
solves (\ref{1:pde2}) for some Lagrange multiplier $\omega \in \R$. Based on these observations, to seek the standing waves of (\ref{new1}) we need only to get solutions of (\ref{1:pde2}), and this can be done by solving the following constrained minimization problem associated with GP energy (\ref{E})
\begin{equation}\label{energy}
d_a(q):=\inf_{\{u\in\mathcal{H},\int_{\mathbb{R}^2} u^2dx=1\}}
E_q(u),
\end{equation}
where  $\mathcal{H}$ is defined by
\begin{equation}\label{H}
\mathcal{H}:=\Big\{u\in
H^1(\mathbb{R}^2):\int_{\mathbb{R}^2}V(x)|u(x)|^2dx<\infty\Big\}.
\end{equation}
Here $V(x):\R^2\to \R ^+$ is
locally bounded and  satisfies $V(x)\to\infty$ as $|x|\to\infty$.
Without loss of generality, by adding a suitable constant
we may assume that
 $$\inf_{x\in
\R^2} V(x) = 0\,,$$
and   $\inf_{x\in
\R^2} V(x)$ can be attained. Under this kind of conditions on $V(x)$, the existence of ground states of (\ref{1:pde2}) was first studied by Rabinowitz \cite{Rab} in some general cases.

Throughout this paper, we denote by $\|u\|_2$ the norm of any functions $u\in L^2(\R^2)$ and $C$ denotes a universal constant which may be different from place to place.

The earlier work related to the minimization problem (\ref{energy}) can
be actually tracked back to the papers \cite{l1,l2,R,Stuart1,Stuart,Stuart2} and the
references therein. A simple scaling argument shows that for the
supercritical case, that is  $q>2$,  (\ref{energy})  does not admit  any
minimizer for all $a>0$. But, in the subcritical case (i.e., $0<q<2$),
(\ref{energy}) admits at least one minimizer for any $a>0$, see e.g.,
\cite{ca,l1,l2}. Moreover, some qualitative properties, such as the
uniqueness, concentration and symmetry, of the minimizers of
(\ref{energy}), for any fixed $0<q<2$, were discussed as
$a\rightarrow+\infty$ in \cite{ca,M} and references therein.
However, for the mass critical case (i.e., $q=2$),  from a physical point of view (cf.\cite{Hulet1,Hulet2,Hulet3}),  there exists a critical cold atom number below which BEC occurs, and collapse occurs otherwise. Mathematically, this was proved  very recently in \cite{Bao,GS}. Roughly speaking,  the authors proved in \cite{Bao,GS}
that there exists a constant $a^*$ such that
(\ref{energy}) admits at least one minimizer if and only if $  a<a^*$,
where
$$a^*:=\|Q\|_2^2,$$
and $Q$ is the unique (up to translations) radially symmetric
positive solution of the following scalar field equation \cite{GNN,K,Li}
\begin{equation}\label{eq1.4}
\Delta u-u+u^3=0 \ \text{ in }\ \mathbb{R}^2,\quad \text{where}\quad
u\in H^1({\mathbb{R}^2}).
\end{equation}
Furthermore, if there are numbers $p_i>0$ and a constant $C>0$ such
that the trapping potential $V(x)$ satisfies
\begin{equation}\label{as:v}
V(x) =  h(x) \prod_{i=1}^n |x-x_i|^{p_i} \quad \text{with $C < h(x)
< 1/C$ for all $x\in \R^2$,}
\end{equation}
the authors in \cite{GS} studied  also the concentration and { symmetry
breaking} of minimizers for $(\ref{energy})$, provided that $q=2$ and $a\nearrow a^*$.

Motivated  by the works mentioned above, in this paper we are interested in
addressing the limit behavior of minimizers for $(\ref{energy})$ when
$q\nearrow 2$ and  $a> a^*$. Towards this
purpose, we first note from \cite{W} that the following scalar field
equation
\begin{equation}\label{eq1.6}
 \Delta u-\frac{2}{q}u+\frac{2}{q}u^{q+1}=0, \ \text{ where } \ q\in (0,2] \ \text{ and }\  u\in H^1(\mathbb{R}^2)
\end{equation}
admits, up to translations, a unique  positive solution which is
radially symmetric about the origin. We denote this unique solution  by $\phi_q=\phi
_q(|x|) $, and throughout the paper, we  set
$$a_q^*:=\|\phi_q\|_2^q.$$
 Moreover, by \cite{W} we have  the following
Gagliardo-Nirenberg inequality
\begin{equation}\label{eq1.7}
\int_{\mathbb{R}^2}|u(x)|^{q+2}dx\leq C_q
\Big\{\int_{\mathbb{R}^2}|\nabla
u(x)|^2dx\Big\}^{\frac{q}{2}}\int_{\mathbb{R}^2}|u(x)|^{2}dx,\quad
u\in H^1({\mathbb{R}^2}),
\end{equation}
where the best constant
$C_q=\frac{q+2}{2\|\phi_q\|_2^q}=\frac{q+2}{2a_q^*}$, and the above
equality holds  at $u(x)=\phi_q(|x|)$.

Note that
$$a_q^*\rightarrow a^*\quad \text{as}\quad q\nearrow2\,,$$
Therefore, for any fixed
$a>a^*$ there exists a constant $\sigma>1$, independent of $q>0$,
such that $\frac{a}{a_q^*}>\sigma >1$ as $q\nearrow 2$, which
further implies that
\begin{equation}\label{1:2.8b}
\Big(\frac{a}{a_q^*}\Big)^{\frac{1}{2-q}}\rightarrow+\infty \ \text{
as } \ q\nearrow 2 .\end{equation} In view of the infinity limit in
(\ref{1:2.8b}), the following {\em   main result} of the present
paper shows the {\em concentration} behavior of minimizers for
(\ref{energy}) as $q\nearrow 2$.

\begin{thm}\label{thm1}
For any fixed  $a>a^*$,  assume that $$V\in C^1(\R^2),\quad
\lim_{|x|\to\infty} V(x) = \infty \quad\text{and} \quad\inf_{x\in \R^2}
V(x) =0\,.$$
Let $u_q\in\mathcal{H}$ be a non-negative minimizer of
(\ref{energy}) with $q\in(0,2)$. Then, for each sequence $\{q_k\}$ with $q_k\nearrow2$ as $k\rightarrow\infty$, there exists a
subsequence of $\{q_k\}$, still denoted by $\{q_k\}$, such that
$u_{q_k}$ concentrates at a global minimum point  $y_0$ of $V(x)$ in
the following sense: for each large $k$, $u_{q_k}$ has a unique
global maximum point $\bar{z}_k\in \R^2$, and satisfies
\begin{equation}\label{1:limt}
\lim_{k\rightarrow\infty}\Big(\frac{a}{a^*_{q_k}}\Big)^{-\frac{1}{2-q_k}}u_{q_k}
\Big(\big(\frac{a}{a^*_{q_k}}\big)^{-\frac{1}{2-q_k}}
x+\bar{z}_k\Big)=\frac{1}{\sqrt{e}\|Q\|_2}Q\Big(\frac{|x|}{\sqrt{e}}\Big)\quad \text{in}\quad H^1(\R^2)\,,
\end{equation} where
$\bar{z}_k\rightarrow y_0$ as $k\rightarrow\infty$.
\end{thm}

Theorem \ref{thm1} gives a detailed description of the behavior of the  minimizers of (\ref{energy}) as $q$ approaches the critical exponent $2$ from
below. Roughly speaking, Theorem \ref{thm1} shows that a minimizer of
(\ref{energy}) behaves like
$$u_{q_k}(x)\approx \Big(\frac{a}{a^*_{q_k}}\Big)^{\frac{1}{2-q_k}}\frac{1}{\sqrt{e}\|Q\|_2}
Q\Big(\frac{(\frac{a}{a^*_{q_k}})^{\frac{1}{2-q_k}}(x-\bar
z_k)}{\sqrt{e}}\Big)\quad \text{as}\ \ q_k\nearrow 2\,.$$  The proof
of Theorem \ref{thm1} is based on precise energy estimates of the GP
energy $d_a(q)$. In fact, we prove in Section 2 [Lemma 2.2] that
\[
 d_a(q)   \thickapprox -\frac{2-q}{2}\Big(\frac{q}{2}\Big)^{\frac{q}{2-q}}
\Big(\frac{a}{a_q^*}\Big)^{\frac{2}{2-q}}  \quad \text{as}\quad
q\nearrow 2\,,
\]
and therefore $ d_a(q) \rightarrow -\infty$ as $q\nearrow 2$ in view
of (\ref{1:2.8b}). As a byproduct of the proof of  Theorem \ref{thm1}, we
shall be able to provide in Lemma \ref{le2.2} the refined
information (compared with those obtained in \cite{ca}) on the minimum energy $\tilde{d}_a(q)$ as well as its
minimizers, where $\tilde{d}_a(q)$ is defined by
\begin{equation*}
\tilde{d}_a(q)=\inf_{\{u\in H^1(\mathbb{R}^2),\int_{\mathbb{R}^2}
u^2dx=1\}} \tilde{E}_q(u)\,,
\end{equation*}
and
\begin{equation}\label{1:eq2.2}
\tilde{E}_q(u):=\int_{\mathbb{R}^2}|\nabla
u(x)|^2dx-\frac{2a}{q+2}\int_{\mathbb{R}^2}|u(x)|^{q+2}dx, \quad
u\in H^1(\mathbb{R}^2)\,.
\end{equation}

Furthermore, we want to show that the concentration point $y_0$ in Theorem
\ref{thm1} is located in the flattest global minimum point of $V(x)$. Towards
this conclusion, we shall assume that the trapping  potential $V(x)$ has
$n\geq 1$ isolated minima, and that $V(x)$ behaves like  in their
vicinity a power of the distance from these points. More precisely,
we shall assume that there exist $n\geq 1$ distinct points $x_i\in
\R^2$ with $V(x_i)=0$, while $V(x)>0$ otherwise. Moreover, there are
numbers $p_i>0$  such that
\begin{equation}\label{as:v}
V(x) =  O( |x-x_i|^{p_i}) \quad \text{near $x_i$,\quad  where
$i=1,2,\cdots ,n$.}
\end{equation}
$\lim_{x\to x_i}
\frac{V(x)}{|x-x_i|^{p_i}}$ exists for all $1\leq i \leq n$.


Let $p = \max \big\{p_1,\cdots ,p_n\big\}$, and let $\lambda_i \in
(0,\infty ]$ be given by
\begin{equation}\label{def:li}
 \lambda_i =    \lim_{x\to x_i} \frac{V(x)}{|x-x_i|^p}  \,.
\end{equation}
Define $\lambda = \min\big\{\lambda_1,\cdots ,\lambda_n\big\}$ and
let
\begin{equation}\label{def:Z}
\mathcal{Z}:=\big\{x_i:\, \lambda_i=\lambda\big\}
\end{equation}
denote the locations of the flattest global minima of $V(x)$. By the
above notations, we have
the following result, which tells us some further information about the concentration point $y_0$ given by in Theorem \ref{thm1}.

\begin{thm}\label{cor}
Under the assumptions of Theorem \ref{thm1} and let $V(x)$ satisfy also the additional condition (\ref{as:v}), then the unique
concentration point $y_0$ obtained in Theorem \ref{thm1} has the properties:
\begin{equation}\label{def:y}
y_0\in\mathcal{Z}\ \text{ and }\
\lim_{k\rightarrow\infty}\big|{\bar{z}_k-y_0}\big|
\Big({\frac{a}{a^*_{q_k}}\Big)^{\frac{1}{2-q_k}}}=0.
\end{equation}
\end{thm}

\begin{rem}
We should mention that if $V(x)$ has some symmetry, for example
$$V(x)=\prod_{i=1}^n|x-x_i|^{p}\quad\text{with}\quad p>0\,,$$ and $x_i$ are arranged on the
vertices of a regular polygon, Theorem \ref{cor} implies the {\em
symmetry breaking} occurring in the minimizers of (\ref{energy}) as
$q\nearrow 2$: there exists  $q_*$ satisfying $0< q_* <2$ such that
for any $q_* < q < 2$, the GP functional (\ref{energy}) has (at
least $n$ different) non-negative minimizers, each of which
concentrates at a specific global minimum point $x_i$. We note that the  symmetry breaking bifurcation for ground states for nonlinear Schr\"odinger or GP equations has been studied in detail in the literature, see, e.g., \cite{J,K08,K11}.
\end{rem}

The results
of the  paper can be extended to general space dimensions $N$
different from 2, if the exponent $q$ in the last term of (\ref{E})
is restricted to the interval  $(0,\frac{4}{N})$, and the limit
$q\nearrow 2$ is replaced by $q\nearrow \frac{4}{N}$. We finally
remark that the concentration phenomena have also been studied
elsewhere in different contexts. For instance, there is a
considerable literature on the concentration phenomena of positive
ground states of the elliptic equation
\begin{equation}
h^2\Delta u(x)-V(x)u(x)+ u^p(x)=0 \quad \mbox{in }\  \R^N
\label{4con:a}
\end{equation}
as $h\to 0^+$, see \cite{BWang,DKW,LW,Wang} and references therein for more
details.

This paper is organized as follows: Section 2 is devoted mainly to
the proof of Theorem \ref{thm2.4} on energy estimates of the
minimizers for (\ref{energy}). We then use Theorem \ref{thm2.4} to
prove Theorem \ref{thm1} in Section 3 by the blow up analysis, and then we prove Theorem \ref{cor} in the end of the section.

\section{Energy Estimates}

The main purpose of this section is to establish Theorem
\ref{thm2.4}, which addresses energy estimates of minimizers for
(\ref{energy}). For any $0<q<2$, let $\phi_q$ be the unique (up to
translations) radially symmetric positive solution of (\ref{eq1.6}).
It then follows
 directly from  Lemma 8.1.2 in \cite{ca} that $\phi_q$ satisfies \begin{equation}\label{eq2.1a}
\int_{\mathbb{R}^2}|\nabla \phi_q(x)|^2dx=\int_{\mathbb{R}^2}|
\phi_q(x)|^2dx=\frac{2}{q+2}\int_{\mathbb{R}^2}| \phi_q(x)|^{q+2}dx.
\end{equation}
Moreover,  one can obtain from \cite{Ber} that there exist positive
constants $\delta$, $C$ and $R_0$, independent of $q>0$, such that
for any $|x|>R_0$,
\begin{equation}\label{eq2.1b}
|\phi_q(x)|+|\nabla \phi_q(x)|\leq C e^{-\delta|x|} \quad
\text{for}\quad q\in[1,2].
\end{equation}
Furthermore, a simple analysis shows that $\phi_q$ satisfies
\begin{equation}\label{eq2.1c}
\phi_q(x)\rightarrow Q(x)\ \text{ strongly in }\ H^1(\mathbb{R}^2)\
\text{ and } \ a_q^*:=\|\phi_q\|_2^q\rightarrow a^*:=\|Q\|_2^2  \
\text{ as }\ q\nearrow 2.
 \end{equation}

We next denote $\tilde{E}_q(u)$ the following energy functional
without the potential
\begin{equation}\label{eq2.2}
\tilde{E}_q(u):=\int_{\mathbb{R}^2}|\nabla
u(x)|^2dx-\frac{2a}{q+2}\int_{\mathbb{R}^2}|u(x)|^{q+2}dx, \quad
u\in H^1(\mathbb{R}^2),
\end{equation}
and consider the associated GP energy
\begin{equation}\label{eq2.3}
\tilde{d}_a(q)=\inf_{\{u\in H^1(\mathbb{R}^2),\int_{\mathbb{R}^2}
u^2dx=1\}} \tilde{E}_q(u).
\end{equation}
It is well known from Chapter 8 in \cite{ca} that if $q\in(0,2)$,
then there exists a unique  (up to translations) positive  minimizer
for $\tilde{d}_a(q)$ at any $a>0$. The following lemma gives the
 refined information on the minimum energy  $\tilde{d}_a(q)$ as well
as its minimizers.

\begin{lem}\label{le2.2} Let $q\in(0,2)$ and $\phi_q$ be the unique radially symmetric positive
 solution of (\ref{eq1.6}). Then,
\begin{equation}\label{eq2.4a}\tilde{d}_a(q)=-\frac{2-q}{2}\Big(\frac{q}{2}\Big)^{\frac{q}{2-q}}
\Big(\frac{a}{a_q^*}\Big)^{\frac{2}{2-q}},\end{equation} and the
unique (up to translations) positive minimizer of $\tilde{d}_a(q)$
must be of the form
\begin{equation}\label{eq2.4}
\tilde{\phi}_q(x)=\frac{\tau_q}{\|\phi_q\|_2}\phi_q(\tau_qx), \quad
\text{where}\quad \tau_q
=\Big(\frac{qa}{2a_q^*}\Big)^{\frac{1}{2-q}}.
\end{equation}
\end{lem}

\noindent{\bf Proof.} By using the Gagliardo-Nirenberg inequality
(\ref{eq1.7}), it follows from (\ref{eq2.2}) that
\begin{equation*}
\tilde{E}_q(u)\geq \int_{\mathbb{R}^2}|\nabla
u(x)|^2dx-\frac{a}{a_q^*}\Big(\int_{\mathbb{R}^2}|\nabla
u(x)|^2dx\Big)^\frac{q}{2},   \text{\ for any \ } u\in H^1(\mathbb{R}^2)\text{ and }
\int_{\mathbb{R}^2} u^2dx=1\,.
\end{equation*}
Let
\begin{equation}\label{eq2.6}
g(s)=s-\frac{a}{a_q^*}s^\frac{q}{2}\quad \text{for}\quad s\in[0,\infty).
\end{equation}
We know that $g(s)$ attains its minimum at $s=\left(\frac{qa}{2a_q^*}\right)^{\frac{2}{2-q}}$, i.e. $s=
\tau_q^2$, which
then implies that
$$\tilde{E}_q(u)\geq g(\tau_q^2) =-\frac{2-q}{2}\Big(\frac{q}{2}\Big)^{\frac{q}{2-q}}\Big(\frac{a}{a_q^*}\Big)^{\frac{2}{2-q}}.$$
This yields that
\begin{equation}\label{eq2.8}
\tilde{d}_a(q)\geq
g(\tau_q^2)=-\frac{2-q}{2}\Big(\frac{q}{2}\Big)^{\frac{q}{2-q}}
\Big(\frac{a}{a_q^*}\Big)^{\frac{2}{2-q}}.
\end{equation}

On the other hand, we introduce the following trial function
$$\psi_q^t(x)=\frac{t}{\|\phi_q\|_2}\phi_q(tx)\quad \text{for}\quad t\in(0,\infty)\,,$$ and $\int_{\mathbb{R}^2}|\psi_q^t|^2dx\equiv
1$ for all $t\in(0,+\infty)$. We then obtain from (\ref{eq2.1a}) that
$$\int_{\mathbb{R}^2}|\nabla \psi_q^t|^2dx=\frac{t^2}{\|\phi_q\|_2^2}\int_{\mathbb{R}^2}|\nabla \phi_q|^2dx=t^2,$$
and
$$\int_{\mathbb{R}^2}|\psi^t_q|^{q+2}dx=\frac{t^q}{\|\phi_q\|_2^{q+2}}\int_{\mathbb{R}^2}| \phi_q|^{q+2}dx=\frac{q+2}{2a_q^*}t^q.$$
Hence
\begin{equation*}
\tilde{d}_a(q)\leq
\tilde{E}_q(\psi_q^t)=t^2-\frac{a}{a_q^*}t^q=g(t^2),\quad \text{ for
any } t\in(0,\infty)\,,
\end{equation*} where $g(\cdot)$ is given by (\ref{eq2.6}).
Thus, we may take $t=\tau_q$, that is,
\begin{equation*}
 \tilde{d}_a(q)\leq g(\tau_q^2)\,,
\end{equation*}
this and (\ref{eq2.8}) then imply the estimate (\ref{eq2.4a}).
Moreover, $ \tilde{d}_a(q)$ is attained at
$\tilde{\phi}_q(x)=\frac{\tau_q}{\|\phi_q\|_2}\phi_q(\tau_qx)$, and
the proof is therefore done in view of the uniqueness (cf.  Chapter
8 in \cite{ca}) of   positive minimizers for  $ \tilde{d}_a(q)$.
\qed\\

\begin{rem}\label{re2.1}
For any fixed
$a>a^*$, since $a_q^*\rightarrow a^*$ as $q\nearrow2$,  there exists a constant $\sigma>1$, independent of $q>0$,
such that $\frac{a}{a_q^*}>\sigma >1$ as $q$ is sufficiently close
to $2^-$. Therefore, we further have
\begin{equation}\label{eq2.8b}
 \tau_q=\Big(\frac{qa}{2a_q^*}\Big)^{\frac{1}{2-q}}\rightarrow+\infty \ \text{ and }\ \tilde{d}_a(q)\rightarrow-\infty\ \text{ as
} \ q\nearrow 2 .\end{equation}
\end{rem}

By applying Lemma \ref{le2.2}, we are able to establish the
following estimates.

\begin{lem}\label{le2.3}
Let $a>a^*$ be fixed, and suppose that
$$V(x)\in L^\infty_{\rm
loc}(\R^2)\,,\quad\lim_{|x|\to\infty} V(x) = \infty\quad\text{and}\quad\inf_{x\in
\R^2} V(x) =0\,.$$ Then,
\begin{equation}\label{eq2.9}
d_a(q)-\tilde{d}_a(q)\rightarrow 0 \quad \text{as}\quad
q\nearrow2\,,
\end{equation}
and
\begin{equation}\label{eeq2.4}
\int_{\mathbb{R}^2}V(x)|u_q(x)|^2dx\rightarrow0 \quad\text{as}\quad
q\nearrow 2\,,
\end{equation}
where $u_q(x)$ is a positive minimizer of (\ref{energy}).
\end{lem}

\noindent{\bf Proof.} By the definitions of $\tilde{d}_a(q)$ and
$d_a(q)$, it is easy to observe that
\begin{equation}\label{eeq2.2}
d_a(q)-\tilde{d}_a(q)\geq 0.
\end{equation}
We next choose  a suitable trial function to estimate the upper bound of $d_a(q)-\tilde{d}_a(q)$. 
For  $R>0$ is fixed, let $\varphi_R(x)\in C_0^\infty(\mathbb{R}^N)$
be a cut-off function
 such that $\varphi_R(x)\equiv1$ if
$x \in B_R(0)$, $\varphi_R(x)\equiv0$ if $x\in B_{2R}^c(0)$, and
$0\leq \varphi_R(x)\leq 1$,
 $|\nabla \varphi(x)|\leq
\frac{C_0}{R}$ for any $x\in B_{2R}(0)\setminus B_R(0)$. Set
\begin{equation}\label{eq2.10}
w_{R,q}(x)=A_{R,q}\tilde{w}_{R,q}(x)=A_{R,q}\varphi_R(x-x_0)\tilde{\phi}_q(x-x_0)\quad
\text{with}\quad x_0\in \mathbb{R}^2\,,
\end{equation}
where $\tilde{\phi}_q(x)$  defined in (\ref{eq2.4}) is the unique
(up to translations) positive minimizer of $\tilde{d}_a(q)$, and
$A_{R,q}>0$ is chosen so that $\|w_{R,q}\|_2^2=1$. It is easy to
calculate that
$$1\leq A_{R,q}^2=\frac{\|\phi_q\|^2_2}{\int_{\mathbb{R}^2}\varphi_R^2(\frac{x}{\tau_q})|\phi_q(x)|^2dx}<\frac{\|\phi_q\|^2_2}{\int_{B_{ R\tau_q}}|\phi_q(x)|^2dx},$$
where $\tau_q>0$ is as in (\ref{eq2.8b}). Since
$\tau_q\rightarrow\infty$ as $q\nearrow2$ and $\phi_q(x)$ decays
 exponentially  as $|x|\to\infty$,
we then have
\begin{equation*}\label{eq2.11}
0\leq A_{R,q}^2-1\leq
\frac{\int_{B^c_{R\tau_q}}|\phi_q(x)|^2dx}{\int_{B_{R\tau_q}}|\phi_q(x)|^2dx}\leq
CR\tau_qe^{-2\delta R\tau_q} \leq Ce^{-\delta R\tau_q}  \quad
\text{as}\quad q\nearrow2,
\end{equation*}
where $\delta >0$ is as in (\ref{eq2.1b}). It hence follows from the
above that
\begin{equation}\label{eeq2.3}
1\leq A_{R,q}^{q+2} \leq (1+Ce^{-\delta
R\tau_q})^{\frac{q+2}{2}}\leq 1+4Ce^{-\delta R\tau_q}.
\end{equation}
In the following, one could take a special value of $R$,  for instance $R = 1$.

Direct calculations show that
\begin{equation}\label{eq2.12}
\begin{split}
&\ \ \ \Big|    \int_{\mathbb{R}^2}|\nabla
\tilde{\phi}_q(x)|^2dx-\int_{\mathbb{R}^2}|\nabla
\tilde{w}_{R,q}(x)|^2dx\Big|\\
&= \Big|\int_{\mathbb{R}^2}|\nabla
\tilde{\phi}_q|^2dx-\int_{\mathbb{R}^2}|\nabla
[\varphi_R(x-x_0)\tilde{\phi}_q(x-x_0)]|^2dx\Big| \\
&=\Big|\int_{\mathbb{R}^2}|\nabla
\tilde{\phi}_q|^2dx-\int_{\mathbb{R}^2}\Big(|\nabla
\varphi_R|^2|\tilde{\phi}_q|^2+|\varphi_R|^2|\nabla \tilde{\phi}_q|^2+2\nabla \varphi_R\varphi_R\nabla \phi_q\phi_q\Big)dx\Big| \\
&\leq
\frac{C}{R^2}\int_{B_R^c}|\tilde{\phi}_q(x)|^2dx+\int_{B_R^c}|\nabla
\tilde{\phi}_q(x)|^2dx+\frac{2C}{R}\int_{B_R^c}|\nabla
\phi_q\|\phi_q|dx.
\end{split}
\end{equation}
Using (\ref{eq2.1b}), we obtain that
\begin{equation}\label{eq2.13}
\begin{split}
&\frac{C}{R^2}\int_{B_R^c}|\tilde{\phi}_q(x)|^2dx=\frac{C}{R^2\|\phi_q\|_2^2}\int_{B_R^c}\tau_q^2|\phi_q(\tau_qx)|^2dx\\
&\leq
\frac{C}{R^2}\int_{B_{R\tau_q}^c}|\phi_q|^2dx<\frac{CR\tau_q}{R^2}e^{-2\delta
R\tau_q}\leq Ce^{-\delta R\tau_q}.
\end{split}
\end{equation}
Similarly,
\begin{equation}\label{eq2.14}
\int_{B_R^c}|\nabla
\tilde{\phi}_q(x)|^2dx=\frac{\tau_q^2}{\|\phi_q\|_2^2}\int_{B_{R\tau_q}^c}|\nabla
\phi_q(x)|^2dx\leq CR\tau_q^3e^{-2\delta R\tau_q}\leq Ce^{-\delta
R\tau_q},
\end{equation}
and
\begin{equation}\label{eq2.15}
\frac{2C}{R}\int_{B_R^c}|\nabla \phi_q\|\phi_q|dx\leq Ce^{-\delta
R\tau_q}.
\end{equation}
It then follows from (\ref{eq2.12})-(\ref{eq2.15}) that
\begin{equation}\label{eq2.16}
 \left|    \int_{\mathbb{R}^2}|\nabla
\tilde{\phi}_q(x)|^2dx-\int_{\mathbb{R}^2}|\nabla
\tilde{w}_{R,q}(x)|^2dx\right|\leq Ce^{-\delta R\tau_q}  \quad
\text{as}\quad q\nearrow2.
\end{equation}
One can also calculate that
\begin{equation}\label{eq2.17}
\left|    \int_{\mathbb{R}^2}|
\tilde{\phi}_q(x)|^{q+2}dx-\int_{\mathbb{R}^2}|
\tilde{w}_{R,q}(x)|^{q+2}dx\right|\leq
\int_{B_R^c}|\tilde{\phi}|^{q+2}dx\leq Ce^{-\delta R\tau_q}.
\end{equation}
Moreover, we have
\[
\int_{R^2}V(x) |w_{R,q}(x)|^2   dx=\frac{A_{R,q}^2}{\|\phi_q\|^2_2}\int V\big(\frac{x}{\tau _q}+x_0\big)\varphi ^2_R\big(\frac{x}{\tau _q}\big)\phi ^2_q (x)dx,
\]
which implies that
\begin{equation*}
\lim_{q\nearrow 2} \int_{R^2}V(x) |w_{R,q}(x)|^2   dx = V(x_0)
\end{equation*}
holds for almost every $x_0\in \R^2$.
Therefore, we choose $x_0\in \mathbb{R}^2$ such that $V(x_0)=0$, and it follows
 from the above estimates that
\begin{equation}\label{2:da}
\begin{split}
&0\leq d_a(q)-\tilde{d}_a(q)\leq E_q(w_{R,q}(x))-\tilde{d}_a(q)\\
&=E_q(A_{R,q}\tilde{w}_{R,q}(x))-\tilde{E}_q(\tilde{\phi}_q(x))\\
&=\Big(E_q(A_{R,q}\tilde{w}_{R,q}(x))-\tilde{E}_q(\tilde{w}_{R,q}(x))\Big)+
\tilde{E}_q(\tilde{w}_{R,q}(x))-\tilde{E}_q(\tilde{\phi}_q(x))\\
&\leq(A_{R,q}^2-1)\int_{\mathbb{R}^2}|\nabla
\tilde{w}_{R,q}|^2dx+\frac{2a}{q+2}(A_{R,q}^{q+2}-1)\int_{\mathbb{R}^2}|
\tilde{w}_{R,q}|^{q+2}dx\\
&\quad +\int_{\mathbb{R}^2}V(x)|w_{R,q}(x)|^2dx  +\left|    \int_{\mathbb{R}^2}|\nabla
\tilde{\phi}_q|^2dx-\int_{\mathbb{R}^2}|\nabla
\tilde{w}_{R,q}(x)|^2dx\right|\\
&\quad +\frac{2a}{q+2}\left|
\int_{\mathbb{R}^2}| \tilde{\phi}_q|^{q+2}dx-\int_{\mathbb{R}^2}|
\tilde{w}_{R,q}(x)|^{q+2}dx\right|\\
&\leq Ce^{-\delta
R\tau_q}+\int_{\mathbb{R}^2}V(x)|w_{R,q}(x)|^2dx\rightarrow0\quad
\text{as}\quad q\nearrow2\,,
\end{split}
\end{equation}
which then implies (\ref{eq2.9}).  By applying the estimate
$$\int_{\mathbb{R}^2}V(x)|u_q(x)|^2dx=d_a(q)-\tilde{E}_q(u_q(x))\leq d_a(q)-\tilde{d}_a(q)\,,$$
we finally conclude (\ref{eeq2.4}) in view of (\ref{eq2.9}). \qed\\

Based on  Lemmas \ref{le2.2} and   \ref{le2.3}, we can establish  the
following delicate estimates.

\begin{thm}\label{thm2.4}
Under the assumptions of Lemma \ref{le2.3},   there exist two
positive constants $C_1$ and $C_2$, independent of $q$, such that
\begin{equation}\label{eq2.20}
\begin{split}C_1\big(\frac{a}{a_q^*}\big)^{\frac{2}{2-q}}&\leq
\int_{\mathbb{R}^2}|\nabla u_q|^2dx     \leq
C_2\big(\frac{a}{a_q^*}\big)^{\frac{2}{2-q}}\quad \text{as}\quad
q\nearrow 2,
\\
C_1\big(\frac{a}{a_q^*}\big)^{\frac{2}{2-q}}&\leq\int_{\mathbb{R}^2}|
u_q|^{q+2}dx\leq C_2\big(\frac{a}{a_q^*}\big)^{\frac{2}{2-q}}\quad
\text{as}\quad q\nearrow 2.\end{split}
\end{equation}
\end{thm}

\noindent{\bf Proof.} By Remark \ref{re2.1} and Lemma \ref{le2.3},
we have $d_a(q)\rightarrow-\infty$ as $q\nearrow2$, and also
\begin{equation}\label{eeq2.5}
\int_{\mathbb{R}^2}|\nabla
u_q|^2dx<\frac{2a}{q+2}\int_{\mathbb{R}^2}| u_q|^{q+2}dx.
\end{equation}
This estimate and the Gagliardo-Nirenberg inequality (\ref{eq1.7})
yield that
$$\frac{2a_q^*}{q+2}\int_{\mathbb{R}^2}| u_q|^{q+2}dx\leq\Big(\int_{\mathbb{R}^2}|\nabla u_q|^2dx\Big)^\frac{q}{2}<\Big(\frac{2a}{q+2}\int_{\mathbb{R}^2}| u_q|^{q+2}dx\Big)^\frac{q}{2},$$
which then implies that
\begin{equation*}
\int_{\mathbb{R}^2}| u_q|^{q+2}dx<\frac{q+2}{2a}
\Big(\frac{a}{a_q^*}\Big)^\frac{2}{2-q}\leq\frac{2}{a}\Big(\frac{a}{a_q^*}\Big)^\frac{2}{2-q}.
\end{equation*}
This establishes the  upper estimates of (\ref{eq2.20}) in view of
(\ref{eeq2.5}).

We address the  lower estimates of (\ref{eq2.20}) as follows. The proof
of Lemma \ref{le2.2} implies that
$$\tilde{d}_a(q)=\tilde{E}_q(\tilde{\phi}_q)=g(s_0),\quad s_0=\tau_q^2=\Big(\frac{q}{2}\Big)^{\frac{2}{2-q}}\Big(\frac{a}{a_q^*}\Big)^{\frac{2}{2-q}},$$
where $g(\cdot )$ is defined as in (\ref{eq2.6}). Since   $g(s)$ is strictly
decreasing in $s\in[0,s_0]$, it follows that for any
$\alpha\in(0,1)$, $$g(s_0)<g(\alpha s_0)<0\quad \text{and}\quad
\gamma_\alpha:=\alpha(-\ln\alpha+1)\in (0,1)\,.$$ Moreover,   direct
calculations show that
\begin{equation*}
0\leq\lim_{q\nearrow 2}\frac{g(\alpha s_0)}{g(s_0)}=\lim_{q\nearrow
2}\frac{\alpha s_0-\frac{a}{a_q^*}(\alpha
s_0)^\frac{q}{2}}{s_0-\frac{a}{a_q^*}s_0^\frac{q}{2}}=
\lim_{q\nearrow2}\frac{2\alpha^\frac{q}{2}-q\alpha}{2-q}=\gamma_\alpha<1,
\end{equation*}
which hence implies that for any $\alpha\in (0,1)$,
\begin{equation}\label{eq2.24}
0>g(\alpha
s_0)>\frac{1+\gamma_\alpha}{2}g(s_0)=\frac{1+\gamma_\alpha}{2}\tilde{d}_a(q)\quad
\text{as}\quad q\nearrow 2.
\end{equation}

We now claim that  for any fixed $0<\alpha<1$, there holds
\begin{equation}\label{eq2.22}
\int_{\mathbb{R}^2}|\nabla u_q|^2dx>\alpha s_0\quad {as}\quad
q\nearrow 2.
\end{equation}
Indeed, if (\ref{eq2.22}) is false, then there exists
$\alpha_0\in(0,1)$, as well as a subsequence of $\{q\}$, still
denoted by $\{q\}$,  such that
$$s_1:=\int_{\mathbb{R}^2}|\nabla u_{q}|^2dx<\alpha_0 s_0\quad \text{as}\quad q\nearrow 2\,.$$
Consequently,
\begin{equation}\label{eq2.25}
d_a(q)= E_{q}(u_q)\geq \int_{\mathbb{R}^2}|\nabla
u_{q}|^2dx-\frac{a}{a^*_{q}}\Big(\int_{\mathbb{R}^2}|\nabla
u_{q}|^2dx\Big)^\frac{q}{2}=g(s_1)>g(\alpha_0 s_0).
\end{equation}
Applying (\ref{eq2.24}), (\ref{eq2.25}) and Lemma \ref{le2.3}, we
then have
\begin{equation*}
\frac{1+\gamma_{\alpha_0}}{2}\tilde{d}_a(q)\leq
d_a(q)<\tilde{d}_a(q)+1,
\end{equation*}
equivalently,
$$\frac{1-\gamma_{\alpha_0}}{2}\tilde{d}_a(q)>-1.$$
This contradicts  the fact that $\tilde{d}_a(q)\rightarrow-\infty$
as $q\nearrow2$. Hence, (\ref{eq2.22}) holds.

Therefore, we obtain the  lower estimates of (\ref{eq2.20}) by applying
(\ref{eeq2.5}) and (\ref{eq2.22}), and the lemma is proved. \qed

\section{Concentration and  Symmetry Breaking}

This section is devoted to  proving Theorem  \ref{thm1} and
Theorem \ref{cor} on the concentration and symmetry breaking of
minimizers for (\ref{energy}) as $q\nearrow 2$, where $a>a^*$ is
fixed. Towards this purpose, we always denote $u_q(x)$ to be a
non-negative minimizer of (\ref{energy}). Set
\begin{equation}\label{eq2.26}
\varepsilon_q :=\varepsilon (q)=\Big(\frac{a}{a_q^*}\Big)^{-\frac{
1}{2-q}}>0\,,
\end{equation}
then $\varepsilon_q \to 0$  by Remark \ref{re2.1}. Define the $L^2(\mathbb{R}^2)$-normalized function
\begin{equation*}
\tilde{w}_q(x):=\varepsilon_q  u_q (\varepsilon_q  x )\,.
\end{equation*}
It then follows from Theorem \ref{thm2.4} that there exist two
positive constants $C_1$ and $C_2$, independent of $q$, such that
\begin{equation}\label{eq2.27}\begin{split}
C_1&\leq \int_{\mathbb{R}^2}|\nabla \tilde{w}_q|^2dx\leq C_2\quad
\text{as}\quad q\nearrow 2\,,\\
 C_1&\leq \int_{\mathbb{R}^2}| \tilde{w}_q|^{q+2}dx\leq C_2\quad
\text{as}\quad q\nearrow 2\,.
\end{split}\end{equation}
We now claim that there exist a sequence $\{y_{\eps_q}\}$,
$R_0>0$ and $\eta>0$ such that
\begin{equation}\label{eq2.28}
\liminf_{\varepsilon_q\rightarrow0}\int_{B_{R_0}(y_{\eps_q})}|\tilde{w}_q|^2dx\geq\eta>0.
\end{equation}
In fact, if (\ref{eq2.28}) is false. Then for any
$R>0$, there exists a sequence $\{\tilde{w}_{q_k}\}$, where
$q_k\nearrow 2$ as $k\to\infty$, such that
\begin{equation*}
\lim_{k\rightarrow\infty}\sup_{y\in
\mathbb{R}^2}\int_{B_{R}(y)}|\tilde{w}_{q_k}|^2dx=0.
\end{equation*}
By Lemma I.1 in \cite{l2} or Theorem 8.10 in \cite{Lieb}, we then
deduce from the above that $\tilde{w}_{q_k}\xrightarrow{k}0$ in
$L^p(\mathbb{R}^2)$ for any $2<p<\infty$. This however contradicts
(\ref{eq2.27}), and the claim is therefore
established.

For the sequence $\{y_{\eps_q}\}$ given by (\ref{eq2.28}), set
\begin{equation}\label{eq2.29}
 w_q(x)=\tilde{w}_q(x+y_{\eps_q} )=\varepsilon_q
u_q(\varepsilon_q x+\varepsilon_q y_{\eps_q}).
\end{equation}
Then  (\ref{eq2.27}) implies that $w_q(x)$ is uniformly
bounded in $H^1(\mathbb{R}^2)$ as $q\nearrow 2$, and the estimate
(\ref{eq2.28}) leads to
\begin{equation}\label{eq2.30}
\liminf_{\varepsilon_q\rightarrow0}\int_{B_{R_0}(0)}|w_q|^2dx\geq\eta>0,
\end{equation}
which therefore implies that $w_q$ cannot vanish as $q\nearrow 2$.

\begin{lem}\label{3:lem1}
Assume $V(x)\in C^1(\R^2)$ satisfies $\lim_{|x|\to\infty} V(x)
= \infty$ and $\inf_{x\in \R^2} V(x) =0$. Then $\{\varepsilon_q
y_{\eps_q}\}$ is uniformly bounded  as $q\nearrow 2$.
Moreover, for any sequence $\{q_k\}$ with $q_k\xrightarrow{k}2$, there exists a
subsequence, still denoted by $\{q_k\}$, such that
$z_k:=\varepsilon_k y_{\varepsilon_k}\xrightarrow{k} y_0$, where $\varepsilon_k:=\varepsilon_{q_k}$ is given by
(\ref{eq2.26}), and $y_0\in \R^2$ is a global minimum point of
$V(x)$, i.e. $V(y_0)=0$.\end{lem}

\noindent{\bf Proof.} It follows from  (\ref{eeq2.4}) and
(\ref{eq2.29}) that
\begin{equation}\label{eq2.31}
\int_{\mathbb{R}^2}V(x)|u_q(x)|^2dx=\int_{\mathbb{R}^2}V(\varepsilon_q
x+\varepsilon_q y_{\eps_q})|w_q(x)|^2dx\rightarrow 0 \quad
\text{as}\quad q\nearrow 2\,.
\end{equation}
Suppose \{$\varepsilon_q y_{\eps_q}$\} is unbounded as $q\nearrow2$, i.e.
$\varepsilon_q\rightarrow0$. Then there exists a subsequence, denoted
by $\{q_n\}$ with $q_n\nearrow2$ as $n\rightarrow\infty$, such that
$$\varepsilon_n:=\varepsilon_{q_n}\xrightarrow{}0 \quad \text{and}\quad \varepsilon_n
y_{\varepsilon_n}\xrightarrow{}\infty \quad \text{as}\quad
n\rightarrow\infty.$$  By the assumptions on $V$, there exists
$C_0>0$ such that $V(x)>C_0$ if $|x|$ is large sufficiently. We then
derive from (\ref{eq2.30}) and Fatou's Lemma that
\begin{equation*}
\lim_{n\rightarrow\infty}\inf \int_{\mathbb{R}^2}V(\varepsilon_n
x+\varepsilon_n y_{\varepsilon_n})|w_{q_n}(x)|^2dx\geq
\int_{\mathbb{R}^2}\liminf_{n\rightarrow\infty}V(\varepsilon_n
x+\varepsilon_n y_{\varepsilon_n})|w_{q_n}(x)|^2dx\geq
{\eta}C_0 >0\,,
\end{equation*}
which however contradicts (\ref{eq2.31}). Thus, $\{\varepsilon_q
y_{\eps_q}\}$ is uniformly bounded for  $q\nearrow2$.
Moreover, for any sequence $\{q_k\}$ with $q_k\xrightarrow{k}2$, there exists a convergent
subsequence, still denoted by $\{q_k\}$, such that
$z_k:=\varepsilon_k y_{\varepsilon_k}\xrightarrow{k} y_0$  for some point $y_0\in \mathbb{R}^2$.

Finally, using  (\ref{eq2.30}) and Fatou's Lemma again, we know that
\begin{equation*}
\lim_{k\rightarrow\infty}\inf \int_{\mathbb{R}^2}V(\varepsilon_k
x+\varepsilon_k y_{\varepsilon_k})|w_{q_k}(x)|^2dx\geq
V(y_0)\int_{B_{R_0}(0)}\lim_{k\rightarrow\infty}|w_{q_k}(x)|^2dx\geq
V(y_0){\eta}\,,
\end{equation*}
which and (\ref{eq2.31}) imply that $V(y_0)=0$, and the lemma  is
therefore proved.\qed\\

Since $u_q$ is a  minimizer of (\ref{energy}), it
satisfies the Euler-Lagrange equation
\begin{equation}\label{eq2.32}
-\Delta u_q(x)+V(x)u_q(x)=\mu_q u_q(x)+a u_q^{q+1}(x)\ \text{ in }\
\mathbb{R}^2,
\end{equation}
where $\mu_q\in \mathbb{R}$ is a  Lagrange multiplier and
satisfies
\begin{equation*}
\mu_q=d_a(q)-\frac{qa}{q+2}\int_{\mathbb{R}^2}|u_q|^{q+2}dx\,.
\end{equation*}
It then follows from Lemma \ref{le2.3} and (\ref{eq2.20}) that there
exist two positive constants $C_1$ and $C_2$, independent of $q$,
such that
$$-C_2<\mu_q \varepsilon_q^2<-C_1\quad \text{as}\quad q\nearrow2\,.$$
By (\ref{eq2.26}) and (\ref{eq2.32}), $w_q(x)$ defined in
(\ref{eq2.29}) satisfies
\begin{equation}\label{eq2.33}
-\Delta w_q(x)+\varepsilon_q^2V(\varepsilon_q x+\varepsilon_q
y_{\eps_q})w_q(x)=\varepsilon_q^2\mu_q w_q(x)+a_q^* w_q^{q+1}(x)\
\text { in }\  \mathbb{R}^2.
\end{equation}
Therefore, by passing to a subsequence if necessary, we can assume
that, for some number $\beta >0$,
$$\mu_{q_k} \varepsilon_k^2 \rightarrow
-\beta^2<0 \quad\text{and}\quad
w_k:=w_{q_k}\rightharpoonup w_0\geq 0 \quad\text{in}\quad H^1(\mathbb{R}^2) \quad\text{as}\quad q_k\nearrow2,$$
for
some $w_0\in H^1(\R^2)$. By passing to the weak limit of
(\ref{eq2.33}), we deduce from Lemma \ref{3:lem1} that the
non-negative function  $w_0$ satisfies
\begin{equation}\label{eq2.34}
-\Delta w(x)=-\beta^2 w(x)+a^* w^{3}(x)\ \text {   in }\
\mathbb{R}^2.
\end{equation}
Furthermore,   we infer from (\ref{eq2.30}) that  $w_0\not\equiv 0$
in $\R^2$, and  the strong maximum principle then yields that
$w_0>0$ in $\R^2$. By a simple rescaling, we thus conclude from the
uniqueness (up to translations) of positive solutions of
(\ref{eq1.4}) that
\begin{equation}\label{eq2.35}
w_0=\frac{\beta}{\|Q\|_2}Q(\beta|x-x_0|) \quad \text{for \, some }\
x_0\in\mathbb{R}^2,
\end{equation}
where $\|w_0\|_2^2=1$. Note that $\|w_k\|_2=1$. Then,  $w_{k}$
converges to $w_0$ strongly in $L^2(\mathbb{R}^2)$ and in fact,
strongly in $L^p(\mathbb{R}^2)$ for any $2\leq p<\infty$ because of
$H^1(\mathbb{R}^2)$ boundedness. Furthermore, since  $w_{k}$ and
$w_0$ satisfy (\ref{eq2.33}) and (\ref{eq2.34}) respectively,
standard elliptic regularity theory gives that  $w_{k}$ converges to
$w_0$
strongly in $H^1(\mathbb{R}^2)$.  \\

\noindent\textbf{Proof of Theorem \ref{thm1}.}  Motivated by
\cite{GS,Wang}, we are now ready to complete the proof of Theorem
\ref{thm1} by the following three steps.
\\

\noindent{\em Step 1: The decay property of $u_k:=u_{q_k}$.} For any
sequence $\{q_k\}$,  let $w_k:=w_{q_k}\ge 0$ be defined by
(\ref{eq2.29}). The above analysis shows that there exists a
subsequence,  still denoted by $\{w_k\}$, satisfying (\ref{eq2.33})
and $w_k\xrightarrow{k} w_0$ strongly in $H^1(\mathbb{R}^2)$ for
some positive function $w_0$. Hence for any $\alpha>2$,
\begin{equation}\label{unif}
\int_{|x|\geq R}|w_k|^\alpha dx\rightarrow0 \quad \text{as}\quad
R\rightarrow\infty \quad \text{uniformly for large }  k.
\end{equation}
Since $\mu_{q_k}<0$, it follows from (\ref{eq2.33}) that $$-\Delta
w_k-c(x)w_k\leq 0\,, \quad\text{where}\quad c(x)=a^*_{q_k}w_k^{q_k}(x)\,.$$ By applying
De Giorgi-Nash-Moser theory (see \cite[Theorem 4.1]{HL}), we thus
have
$$\max_{B_1(\xi)} w_k\leq C\Big(\int_{B_2(\xi)}|w_k|^\alpha
dx\Big)^\frac{1}{\alpha},$$ where $\xi$ is an arbitrary point in
$\mathbb{R}^2$, and $C$ is a constant depending only on the bound of
$\|w_k\|_{L^\alpha(B_2(\xi))}$. We hence deduce from (\ref{unif})
that
\begin{equation}
w_k(x)\rightarrow 0 \quad \text{as}\quad |x|\rightarrow\infty \quad
\text{uniformly in } k.
\end{equation}
Since $w_k$ satisfies (\ref{eq2.33}), one can use  the comparison
principle as in \cite{KW} to compare $w_k$ with
$Ce^{-\frac{\beta}{2}|x|}$, which then shows that there exists a
large constant $R>0$, independent of $k$, such that
\begin{equation}\label{exp}
w_k(x)\leq Ce^{-\frac{\beta}{2}|x|} \quad \text{for} \quad |x|>R
\quad \text{as}\quad k\rightarrow\infty.
\end{equation}
By Lemma \ref{3:lem1}, we therefore obtain from (\ref{exp}) that the
subsequence
$$u_k(x):=u_{q_k}(x)=\frac{1}{\varepsilon_k}w_k(\frac{x-z_k}{\varepsilon_k})$$
 decays uniformly to zero for $x$ outside any fixed neighborhood of
$y_0$ as $k\rightarrow\infty$, where $\varepsilon_k=\varepsilon_{q_k}$, $z_k\in \R^2$ is defined as in Lemma \ref{3:lem1}, and $y_0\in \R^2$ is a global minimum point of $V(x)$.\\

\noindent{\em Step 2: The detailed concentration behavior.} Let
$\bar{z}_k$ be any local maximum point of $u_k$. It then yields from
(\ref{eq2.32})   that $$u_k(\bar{z}_k)\geq
(\frac{-\mu_{q_k}}{a})^\frac{1}{q_k}\geq C\varepsilon_k^{-1}\,.$$ This
estimate and the above decay property thus imply that
$\bar{z}_k\rightarrow y_0$ as $k\rightarrow \infty$.  Set
\begin{equation}\label{3:W}\bar{w}_k=\varepsilon_ku_k(\varepsilon_k x+\bar{z}_k),\end{equation}
so that $\bar{w}_k$ satisfies (\ref{eq2.27}). It then follows from
(\ref{eq2.32}) that
\begin{equation}\label{3:222}
-\Delta \bar{w}_k(x)+\varepsilon_k^2V(\varepsilon_k
x+\bar{z}_k)\bar{w}_k(x)=\varepsilon_k^2\mu_{q_k}
\bar{w}_k(x)+a^*_{q_k} \bar{w}_k^{q_k+1}(x) \ \text { in } \
\mathbb{R}^2.
\end{equation}
The same argument as proving (\ref{eq2.34}) yields that there exists
a subsequence of $\{\bar{w}_k\}$, still denoted by $\{\bar{w}_k\}$,
such that $\bar{w}_k\xrightarrow{k} \bar{w}_0$ in $H^1(\mathbb{R}^2)$
for some nonnegative function $\bar{w}_0\geq 0$, where $\bar w_0$
satisfies (\ref{eq2.34}) for some  constant $\beta >0$. We derive
from (\ref{3:222})  that
\begin{equation}\label{bound}
\bar{w}_k(0)\geq \Big(\frac{-\varepsilon
^2_k\mu_{q_k}}{a^*_{q_k}}\Big)^\frac{1}{q_k}\geq
\Big(\frac{\beta^2}{2a^*}\Big)^\frac{1}{2}\quad \text{as}\quad
k\rightarrow\infty,
\end{equation}
which implies that $\bar{w}_0(0)\geq
(\frac{\beta^2}{2a^*})^\frac{1}{2}$. Thus, the strong maximum
principle yields that $\bar{w}_0(x)>0$ in $\R^2$. Since the $x=0$ is
a critical point of $\bar{w}_k$ for all $k>0$, it is also a critical
point of $\bar{w}_0$.  We therefore conclude from the uniqueness (up
to translations) of positive radial solutions for (\ref{eq1.4}) that
$\bar{w}_0$ is spherically symmetric about the origin, and
\begin{equation}\label{eq2.350}
\bar w_0=\frac{\beta}{\|Q\|_2}Q(\beta|x|) \quad \text{for \, some }\
\beta >0.
\end{equation}

One can deduce from the above that $\bar{w}_k\geq
(\frac{\beta^2}{2a^*})^\frac{1}{2}$ at each local maximum point.
Since $\bar{w}_k$ decays to zero uniformly in $k$ as
$|x|\rightarrow\infty$, all   local maximum points of $\bar{w}_k$
stay in a finite ball in $\mathbb{R}^2$. Since
$\bar{w}_k\xrightarrow{k} \bar{w}_0$ in $C^2_{loc}(\mathbb{R}^2)$
and $x=0$ is the only critical point of $\bar{w}_0$, all local
maximum points must approach the origin and hence stay in a small
ball $B_\epsilon(0)$ as $k\rightarrow\infty$. One can take
$\epsilon$ small enough such that $\bar{w}''_0(r)<0$ for $0\leq
r\leq \epsilon$. It then follows from Lemma 4.2 in \cite{NT} that
for large $k$, $\bar{w}_k$ has no critical points other than the
origin. This  gives the uniqueness of local maximum points
for $\bar{w}_k(x)$, which therefore implies that there exists a unique maximum point $\bar z_k$ for each $\{u_k\}$  and $\{\bar z_k\}$ goes to a  global minimum point of potential $V(x)$ as $k\to\infty$.\\

\noindent{\em Step 3: The exact  value of $\beta$ defined in
(\ref{eq2.350}).} Let $\{q_k\}$, where $q_k\nearrow 2$ as
$k\to\infty$, be the subsequence obtained in Step 2, and denote
$u_k:=u_{q_k}$.
 Recall from Lemma \ref{le2.3} that
$$d_a(q_k)=\tilde{d}_a(q_k)+o(1)= -\frac{2-q_k}{2} \Big(\frac{q_k}{2}\Big)^{\frac{q_k}{2-q_k}}
\varepsilon _k^{-2}+o(1) \quad \text{as}\quad k\to\infty,$$ which
yields that
\begin{equation}\label{eq2.36}
\lim_{k\to\infty}\frac{2}{2-q_k}\varepsilon _k^2
d_a(q_k)=-\lim_{k\to\infty}
\big(\frac{q_k}{2}\big)^{\frac{q_k}{2-q_k}}=-e^{-1}.
\end{equation}
On the other hand,
\begin{equation}\label{eq2.42}
\begin{split}
d_a(q_k)&=\int_{\mathbb{R}^2}|\nabla
u_k|^2dx-\frac{2a}{q_k+2}\int_{\mathbb{R}^2}|u_k|^{q_k+2}dx+\int_{\mathbb{R}^2}V(x)|u_k|^2dx\\
&=\varepsilon _k^{-2}\Big[\int_{\mathbb{R}^2}|\nabla
\bar w_k|^2dx-\frac{2a_{q_k}^*}{q_k+2}\int_{\mathbb{R}^2}|\bar w_k|^{q_k+2}dx\Big]+\int_{\mathbb{R}^2}V(x)|u_k|^2dx\\
&\geq \varepsilon _k^{-2}\Big[\int_{\mathbb{R}^2}|\nabla \bar
w_k|^2dx-\Big(\int_{\mathbb{R}^2}|\nabla \bar
w_k|^2dx\Big)^\frac{q_k}{2}\Big],
\end{split}
\end{equation}
where $\bar w_k:=\bar w_{q_k}$ is as in (\ref{3:W}). Set
$\beta_{q_k}^2:=\int_{\mathbb{R}^2}|\nabla \bar w_k|^2dx$. Since
$\bar w_k(x)\xrightarrow{k} \bar w_0(x)$ strongly in
$H^1(\mathbb{R}^2)$, we have
\begin{equation}\label{3:beta}\lim_{k\to\infty}\beta_{q_k}^2=\|\nabla \bar w_0\|_2^2=\beta ^2,\end{equation}
where (\ref{eq2.1a}) is used. Let $f_k(t)=t-t^\frac{q_k}{2}$, where
$ t\in(0,\infty)$. The simple analysis shows that  $f_k(\cdot)$
attains  its global minimum at the unique point
$t_k:=(\frac{q_k}{2})^\frac{2}{2-q_k}$, and also
$f_k(t_k)=-\frac{2-q_k}{2}(\frac{q_k}{2})^\frac{q_k}{2-q_k}$. We
hence deduce from (\ref{eq2.42}) that
\begin{equation*}
\lim_{k\to\infty}\frac{2}{2-q_k}\varepsilon_k^2d_a(q_k)\geq
\lim_{k\to\infty}\frac{2}{2-q_k}f_k(\beta_k^2)\geq
\lim_{k\to\infty}\frac{2}{2-q_k}f_k(t_k)=-e^{-1},
\end{equation*}
which and (\ref{eq2.36}) lead to the limit
\begin{equation*}
\lim_{k\to\infty}f_k(\beta_k^2)\big/ f_k(t_k)=1.
\end{equation*}
We then obtain that
$$\lim_{k\to\infty}\beta_k^2=\lim_{k\to\infty}t_k=e^{-1},$$
and therefore we have $\beta =e^{-\frac{1}{2}}$ by applying
(\ref{3:beta}), which, together with (\ref{3:W}) and (\ref{eq2.350})
give (\ref{1:limt}). We thus complete the proof of Theorem
\ref{thm1}. \qed\\

Following the proof of Theorem \ref{thm1}, we next address
Theorem \ref{cor} on the local properties of concentration points. Under the
assumption (\ref{as:v}), we first denote
\[\bar{V}_i(x)=V(x)/|x-x_i|^{p_i},\quad \text{where }\ i=1,\cdots,n,\]
so that the limit $\lim_{x\to x_i}\bar{V}_i(x)=\bar{V}_i(x_i)$ is assumed to exist for all $i=1,\cdots ,n$.\\

\noindent\textbf{Proof of Theorem \ref{cor}.} For convenience we still
denote $\{q_k\}$ to be the subsequence obtained in Theorem \ref{thm1}.
Choose a point $x_{i_0}\in \mathcal{Z}$, where $\mathcal{Z}$ is
defined by (\ref{def:Z}), and let
$$w_{R,q_k}(x)=A_{R,q_k}\varphi_R(x-x_{i_0})\tilde{\phi}_{q_k}(x-x_{i_0})$$
be the trial function defined by (\ref{eq2.10}). By (\ref{2:da}), we
know that
\begin{equation}\label{eq2.43}
\begin{split}
d_a(q_k)-\tilde{d}_a(q_k)&\leq
E\big(w_{R,q_k}(x)\big)-\tilde{E}\big(\tilde{\phi}_{q_k}(x-x_{i_0})\big)\\
&\leq
\int_{\mathbb{R}^2}V(x)|w_{R,q_k}(x)|^2dx+Ce^{-\delta R\tau_{q_k}}\\
&\leq\frac{A_{R,q_k}^2}{\tau_{q_k}^p\|\phi_{q_k}\|_2^2}\int_{B_{2R\tau_{q_k}}}\bar{V}_{i_0}
\big(\frac{x}{\tau_{q_k}}+x_{i_0}\big)|x|^p\phi_{q_k}^2(x)dx+Ce^{-\delta
R\tau_{q_k}}\\
&=\frac{A_{R,q_k}^2}{\tau_{q_k}^p\|\phi_{q_k}\|_2^2}\inte\chi_{B_{2R\tau_{q_k}}}(x) \bar{V}_{i_0}
\big(\frac{x}{\tau_{q_k}}+x_{i_0}\big)|x|^p\phi_{q_k}^2(x)dx+Ce^{-\delta
R\tau_{q_k}}
\end{split}
\end{equation}
where $\tau_{q_k}>0$ satisfies
$\tau_{q_k}=(\frac{q_k}{2})^{\frac{1}{2-q_k}}\frac{1}{\varepsilon_k}$ in view
of Lemma \ref{le2.2} and (\ref{eq2.26}), and $\chi_{B_{2R\tau_{q_k}}}$ is the characteristic function of the set $B_{2R\tau_{q_k}}$.
Since $\phi_{q_k}(x)$ decays
exponentially and $\phi_{q_k}\rightarrow Q$ strongly in
$L^2(\mathbb{R}^2)$, then,

$$\chi_{B_{2R\tau_{q_k}}}(x)\bar{V}_{i_0}
(\frac{x}{\tau_{q_k}}+x_{i_0})|x|^p\phi_{q_k}^2(x)\leq \sup_{B_{2R}}\bar{V}_{i_0}
(x+x_{i_0})\cdot C e^{-\delta |x|}\in L^1(\R^2)\,, $$ and
$$\chi_{B_{2R\tau_{q_k}}}(x)\bar{V}_{i_0}
(\frac{x}{\tau_{q_k}}+x_{i_0})|x|^p\phi_{q_k}^2(x)\to \bar{V}_{i_0}
(x_{i_0}) |x|^pQ^2(x)\quad \text{a.e.}\quad \R^2\quad \text{as}\quad k\to\infty\,.$$
Note that $A_{R,q_k}\to 1$ as $q_k\nearrow 2$,
 we thus
obtain from (\ref{eq2.43}) and Lebesgue's dominated convergence theorem that
\begin{equation}\label{eq2.44}
\begin{split}
&\lim_{k\rightarrow
\infty}\frac{d_a(q_k)-\tilde{d}_a(q_k)}{\varepsilon_k^p}\\
&\leq
\lim_{k\rightarrow \infty}\Big(\frac{q_k}{2}\Big)^\frac{-p}{2-q_k}\Big[
\frac{A_{R,q_k}^2}{\|\phi_{q_k}\|_2^2}\inte \chi_{B_{2R\tau_{q_k}}}(x)\bar{V}_{i_0}
(\frac{x}{\tau_{q_k}}+x_{i_0})|x|^p\phi_{q_k}^2(x)dx+C\tau_{q_k}^pe^{-\delta
R\tau_{q_k}}\Big]\\
&=\frac{e^\frac{p}{2}}{\|Q\|_2^2}\lim_{k\rightarrow \infty}\int_{\mathbb{R}^2}\chi_{B_{2R\tau_{q_k}}}(x)\bar{V}_{i_0}
(\frac{x}{\tau_{q_k}}+x_{i_0})|x|^p\phi_{q_k}^2dx\\
&=\frac{\bar{V}_{i_0}(x_{i_0})e^\frac{p}{2}}{\|Q\|_2^2}\int_{\mathbb{R}^2}|x|^pQ^2dx.
\end{split}
\end{equation}

On the other hand, following the proof of Theorem \ref{thm1} we
denote $\bar z_k$ to be the unique global maximum point of $u_k$,
and let $\bar w_k$ be defined as in (\ref{3:W}). Denote also $y_0\in
\R^2$ to be the limit of $\bar z_k$ as $k\to\infty$. Since
$V(y_0)=0$, then there exists an $x_j=y_0$ for some $1\leq j\leq n$. We claim that  $\big\{\frac{\bar{z}_k-x_j}{\varepsilon_k}\big\}$ is
uniformly bounded in $\mathbb{R}^2$. Indeed, if there exists a
subsequence, still denoted by $\{q_k\}$, such that
$\frac{\bar{z}_k-x_j}{\varepsilon_k}\rightarrow\infty$  as
$k\to\infty$, it then follows from Fatou's Lemma that,  for any $C>0$ sufficiently large,
\begin{equation}
\begin{split}
\lim_{k\to\infty}\frac{d_a(q_k)-\tilde{d}_a(q_k)}{\varepsilon
_k^{p_j}}&\geq\lim_{k\to\infty}\int_{\mathbb{R}^2}\bar{V}_j(\varepsilon_k
x+\bar{z}_k)\Big|x+\frac{\bar{z}_k-x_j}{\varepsilon_k}\Big|^{p_j}\bar
w_k^2dx\\
&\geq\int_{\mathbb{R}^2}\lim_{k\to\infty}\bar{V}_j(\varepsilon_k
x+\bar{z}_k)\Big|x+\frac{\bar{z}_k-x_j}{\varepsilon_k}\Big|^{p_j}\bar
w_k^2dx\geq C \bar{V}_j(x_j)\,,
\end{split}
\end{equation}
which however contradicts (\ref{eq2.44}) owing to  $p_j\leq p=\max \big\{p_1,\cdots ,p_n\big\}$, and the claim is therefore
true. Consequently, there exists a subsequence,  still denoted by
$\{q_k\}$, such that
\begin{equation}\label{eq2.46a}
\frac{\bar{z}_k-x_j}{\varepsilon_k}\rightarrow \bar{z}_0\quad
\text{for some }\ \bar{z}_0\in\mathbb{R}^2.\end{equation} Since $Q$
is a radial decreasing function and decays exponentially as
$|x|\to\infty$, we then deduce that
\begin{align}
\lim_{k\to\infty}
\frac{d_a(q_k)-\tilde{d}_a(q_k)}{\varepsilon_k^{p_j}}&
\geq\lim_{k\to\infty}\int_{\mathbb{R}^2}\bar{V}_j(\varepsilon_k
x+\bar{z}_k)\left|x+\frac{\bar{z}_k-x_j}{\varepsilon_k}\right|^{p_j}\bar w_k^2dx\nonumber\\
&\geq \bar{V}_j(x_j)\int_{\mathbb{R}^2} |x+\bar{z}_0|^{p_j}\bar{w}_0^2dx\nonumber\\
&=\frac{\bar{V}_{j}(x_j)e^\frac{p_j}{2}}{\|Q\|_2^2}
\int_{\mathbb{R}^2}\big|x+\frac{\bar{z}_0}{\sqrt{e}}\big|^{p_j}Q^2dx\nonumber\\
&\geq\frac{\bar{V}_{j}(x_j)e^\frac{p_j}{2}}{\|Q\|_2^2}
\int_{\mathbb{R}^2}|x|^{p_j}Q^2dx,\label{eq2.46}
\end{align}
where $\bar{w}_0>0$ is as in (\ref{eq2.350}), and ``=" in the last inequality of (\ref{eq2.46})
holds  if and only if $\bar z_0=(0,0)$.

From (\ref{eq2.44}) and (\ref{eq2.46}), we see that $p_j\geq p$, however, since $p = \max \big\{p_1,\cdots ,p_n\big\}$, we thus have
$p_j=p$. And then, by comparing  (\ref{eq2.44}) with (\ref{eq2.46}) again, we get that
$\bar{V}_j(x_j)\leq \bar{V}_{i_0}(x_{i_0})$. Meanwhile,
$\bar{V}_j(x_j)\geq \bar{V}_{i_0}(x_{i_0})$ always holds for $x_{i_0}\in\mathcal{Z}$. Thus, $\bar{V}_j(x_j)=
\bar{V}_{i_0}(x_{i_0})$, this means that
$x_j=y_0\in\mathcal{Z}$ must be the flattest global minimum point of
$V(x)$. These  further  yield that   (\ref{eq2.46}) is indeed an
equality, therefore $\bar z_0=(0,0)$, which gives (\ref{def:y}).  This
completes the proof of  Theorem \ref{cor}. \qed

{\bf Acknowledgements:} This work was supported by the National Natural Science Foundation of China
 under grants No.11071245, 11241003, 11322014 and  11271360.


\begin{thebibliography}{GNN}

\bibitem{Bao} W. Z. Bao  and Y. Y. Cai, {\em Mathematical theory and numerical methods for Bose-Einstein condensation}, Kinet. Relat. Models {\bf 6} (2013), no. 1, 1--135.

\bibitem{Ber} H. Berestycki  and P. L. Lions, {\em Nonlinear scalar field equations. I. Existence of a ground state}, Arch. Rat. Mech. Anal. {\bf 82} (1983), 313--346.


\bibitem{Hulet1} C. C. Bradley, C. A. Sackett, J. J. Tollett and R. G. Hulet, {\it Evidence of Bose-Einstein condensation in an atomic gas with attractive interactions}, Phys. Rev. Lett. {\bf 75} (1995), 1687--1690. {\it Erratum} Phys. Rev. Lett. {\bf 79} (1997), 1170.

\bibitem{Hulet2} C. C. Bradley, C. A. Sackett and R. G. Hulet, {\it Bose-Einstein condensation of lithium: observation of limited condensate number}, Phys. Rev. Lett. {\bf 78}  (1997), 985--989.

\bibitem{BWang} J. Byeon and Z. Q. Wang, {\em Standing waves with a critical frequency for nonlinear Schr\"{o}dinger equations}, Arch. Ration. Mech. Anal. {\bf 165} (2002), no. 4, 295--316.

\bibitem{ca} T. Cazenave, Semilinear Schr\"{o}dinger Equations, Courant Lecture Notes in Mathematics  Vol. 10, Courant Institute of Mathematical Science/AMS, New York, 2003.






\bibitem{DKW} M. del Pino, M. Kowalczyk and J. C. Wei, {\it Concentration on curves for nonlinear schr\"odinger equations }, Comm. Pure Appl. Math. {\bf 60}, (2007), 113--146.


\bibitem{GNN} B. Gidas, W. M. Ni and L. Nirenberg, {\em Symmetry of positive solutions of nonlinear elliptic equations in $\R^n$}, Mathematical analysis and applications  Part A, Adv. in Math. Suppl. Stud. vol. {\bf 7}, Academic Press, New York, (1981), 369--402.


\bibitem{GS} Y. J. Guo and R. Seiringer, {\em On the mass concentration of Bose-Einstein condensates with attactive interactions},
preprint, (2012).


\bibitem {HL}  Q. Han and F. H. Lin,   Elliptic Partial Differential Equations: Second Edition, Courant Lecture Notes in Mathematics  Vol. 1, Courant Institute of Mathematical Science/AMS, New York, 2011.


\bibitem {J} R. K. Jackson and M. I. Weinstein, {\em Geometric analysis of bifurcation and symmetry breaking in a Gross-Pitaevskii equation}, J. Stat. Phys. {\bf 116} (2004), 881--905.


\bibitem{KW} O. Kavian and F. B. Weissler, {\em Self-similar solutions of the pseudo-conformally invariant nonlinear Schr\"{o}dinger equation},
Michigan Math. J. {\bf 41} (1994), no. 1, 151--173.

\bibitem {K08} E.W. Kirr, P.G. Kevrekidis,  E. Shlizerman and  M.I.  Weinstein,  {\em  Symmetry-breaking bifurcation in nonlinear Schr\"odinger/Gross-Pitaevskii equations}, SIAM J. Math. Anal. {\bf 40} (2008), 566--604.

 \bibitem {K11}  E. W. Kirr, P. G. Kevrekidis and D. E. Pelinovsky, {\em
Symmetry-breaking bifurcation in the nonlinear Schr\"odinger equation with symmetric potentials}, Comm. Math. Phys. {\bf 308} (2011), no. 3, 795--844.


\bibitem{K} M. K. Kwong, {\em Uniqueness of positive solutions of $\Delta u-u+u^p=0$ in $\R^N$}, Arch. Rational Mech. Anal. {\bf 105} (1989), 243--266.



\bibitem{Li} Y. Li and W.-M. Ni, {\em Radial symmetry of positive solutions of nonlinear elliptic equations in $\R^n$}, Comm. Partial Differential Equations {\bf 18} (1993), 1043--1054.

\bibitem{Lieb} E. H. Lieb and M. Loss, Analysis, Graduate Studies in Mathematics Vol. 14. Amer. Math. Soc., Providence, RI, second edition, 2001.



\bibitem {LSY} E. H. Lieb, R. Seiringer and J. Yngvason, {\em Bosons in a trap: A rigorous derivation of the Gross-Pitaevskii energy functional}, Phys. Rev. A {\bf 61}, 043602-1-13 (2000).

\bibitem{l1} P. L. Lions, {\em The concentration-compactness principle in the caclulus of variations. The locally compact case I}, Ann. Inst. H. Poincar\'{e} Anal. Non Lin\'{e}aire. {\bf 1} (1984), 109--145.

\bibitem{l2} P. L. Lions, {\em The concentration-compactness principle in the caclulus of variations. The locally compact case II}, Ann. Inst. H. Poincar\'{e} Anal. Non Lin\'{e}aire. {\bf 1} (1984), 223--283.

\bibitem{LW} G. Z. Lu and J. C. Wei, {\em On nonlinear schr\"odinger equations with totally degenerate potentials}, C. R. Acad. Sci. Paris. {\bf 326} (1998), 691--696.


\bibitem {NT} W.-M. Ni and I. Takagi, {\em On the shape of least-energy solutions to a semilinear Neumann problem}, Comm. Pure Appl. Math. {\bf 44} (1991), 819--851.


\bibitem {M} M. Maeda, {\em On the symmetry of the ground states of nonlinear Schr\"{o}dinger equation with potential}, Adv. Nonlinear Stud. {\bf 10} (2010), no. 4, 895--925.


\bibitem{Rab} P. H.  Rabinowitz, {\em On a class of nonlinear Schr\"odinger equations}, Z. Angew. Math. Phys. {\bf 43} (1992), 270--
291.

\bibitem {RS} M. Reed and B. Simon, Methods of modern mathematical physics. IV. Analysis of operators, Academic Press, New York-London, 1978.

\bibitem {R}  H. A.  Rose and M.I.  Weinstein, {\em On the bound states of the nonlinear Schr\"odinger equation with a linear potential}, Physica D {\bf 30} (1988), 207--218.


\bibitem{Hulet3} C. A. Sackett, H. T. C. Stoof and R. G. Hulet, {\it Growth and collapse of a Bose-Einstein condensate with attractive interactions}, Phys. Rev. Lett. {\bf 80} (1998), 2031--2034.

\bibitem {S} R. Seiringer, {\em Hot topics in cold gases}, XVIth International Congress on Mathematical Physics, World Sci. Publ., Hackensack, NJ, (2010), 231--245.

\bibitem {Stuart1} C. A.  Stuart, {\em Bifurcation for Dirichlet problems without eigenvalues},   Proc. London Math. Soc., {\bf 45} (1982), 169-192.

 \bibitem {Stuart} C. A.  Stuart, {\em Bifurcation from the essential spectrum},   Lecture Notes in Math. {\bf 1017}, Springer, Berlin, 1983.

    \bibitem {Stuart2} C. A.  Stuart, {\em Bifurcation from the essential spectrum for some non-compact non-linearities },   Math. Methods Applied Sci., {\bf 11} (1989), 525-542.



\bibitem {Wang} X. F. Wang, {\em On concentration of positive bound states of nonlinear Schr\"{o}dinger equations}, Comm. Math. Phys. {\bf 153} (1993), no. 2, 229--244.

\bibitem {W} M. I. Weinstein, {\em Nonlinear Schr\"{o}dinger equations and sharp interpolations estimates}, Comm. Math. Phys. {\bf 87} (1983), 567--576.







\end{thebibliography}
\end{document}